\documentclass[12pt,reqno]{amsart}
\usepackage{latexsym, amsmath, amscd, amssymb, amsthm, bm, dsfont}
\usepackage[l >= 1em]{diagrams}
\usepackage[colorlinks=true, linkcolor=blue,urlcolor=blue,
  citecolor=blue]{hyperref}
\usepackage[shortalphabetic]{amsrefs}
\usepackage[T1]{fontenc}
\usepackage{mathptmx}
\usepackage{microtype}

\usepackage[centering, includeheadfoot, hmargin=1.2in, vmargin=0.8in,
  headheight=30.4pt]{geometry}
\newtheorem{lemma}{Lemma}[section]
\newtheorem{theorem}[lemma]{Theorem}
\newtheorem{corollary}[lemma]{Corollary}
\newtheorem{proposition}[lemma]{Proposition}
\newtheorem{conjecture}[lemma]{Conjecture}
\theoremstyle{definition}

\newtheorem{example}[lemma]{Example}
\newcommand{\define}[1]{{\bfseries\itshape #1}}
  
\makeatletter
\renewcommand\subsection{\@startsection{subsection}{2}%
  \z@{.5\linespacing\@plus.7\linespacing}{-.5em}%
  {\normalfont\scshape}}
\makeatother
\newcommand{\relphantom}[1]{\mathrel{\phantom{#1}}}
\newcommand{\kk}{\ensuremath{\Bbbk}}

\newcommand{\NN}{\ensuremath{\mathbb{N}}} 
\newcommand{\PP}{\ensuremath{\mathbb{P}}} 
 
\newcommand{\RR}{\ensuremath{\mathbb{R}}} 
\newcommand{\ZZ}{\ensuremath{\mathbb{Z}}} 
\newcommand{\sF}{\ensuremath{\mathcal{F}}} 
\newcommand{\sO}{\ensuremath{\mathcal{O}}} 

\newcommand{\fB}{\ensuremath{\mathfrak{B}}} 
\newarrow{Equal}{=}{=}{=}{=}{=}
\newcommand{\id}{\ensuremath{\mathds{1}}}
\renewcommand{\geq}{\geqslant}
\renewcommand{\leq}{\leqslant}
\DeclareMathOperator{\Ideal}{I}
\DeclareMathOperator{\Ker}{Ker}
\DeclareMathOperator{\Nef}{Nef}
\DeclareMathOperator{\Pic}{Pic}
\DeclareMathOperator{\pos}{pos}
\DeclareMathOperator{\Proj}{Proj}
\DeclareMathOperator{\Sym}{Sym}
\DeclareMathOperator{\conv}{conv}

\begin{document}

\title[Linear Determinantal Equations]{Linear Determinantal Equations for all
  Projective Schemes}

\author[J.~Sidman]{Jessica Sidman}
\address{Department of Mathematics and Statistics\\ Mount Holyoke College\\
  South Hadley\\ MA\\ 01075\\ USA}
\email{jsidman@mtholyoke.edu}

\author[G.G.~Smith]{Gregory G. Smith} 
\address{Department of Mathematics and Statistics \\ Queen's University \\
  Kingston \\ ON \\ K7L~3N6\\ Canada}
\email{ggsmith@mast.queensu.ca}



\begin{abstract}
  We prove that every projective embedding of a connected scheme determined by
  the complete linear series of a sufficiently ample line bundle is defined by
  the $2 \times 2$ minors of a $1$-generic matrix of linear forms.  Extending
  the work of Eisenbud-Koh-Stillman for integral curves, we also provide
  effective descriptions for such determinantally presented ample line bundles
  on products of projective spaces, Gorenstein toric varieties, and smooth
 varieties.
\end{abstract}

\maketitle

\section{Introduction} 
\label{s:intro}

\noindent
Relating the geometric properties of a variety to the structural features of its
defining equations is a fundamental challenge in algebraic geometry. Describing
generators for the homogeneous ideal associated to a projective scheme is a
basic form of this problem. For a rational normal curve, a Segre variety, or a
quadratic Veronese variety, the homogeneous ideal is conveniently expressed as
the $2$-minors (i.e.{} the determinants of all $2 \times 2$ submatrices) of a
generic Hankel matrix, a generic matrix, or a generic symmetric matrix
respectively.  These determinantal representations lead to a description of the
minimal graded free resolution of the homogeneous ideal of the variety and
equations for higher secant varieties. Mumford's ``somewhat startling
observation'' in \cite{Mumford}*{p.~31} is that a suitable multiple of every
projective embedding is the intersection of a quadratic Veronese variety with a
linear space and, hence, defined by the $2$-minors of a matrix of linear
forms. Exercise~6.10 in \cite{Eisenbud} rephrases this as a ``(vague) principle
that embeddings of varieties by sufficiently positive bundles are often defined
by ideals of $2 \!\times\! 2$ minors''. Our primary goal is to provide a precise
form of this principle.

To be more explicit, consider a scheme $X$ embedded in $\PP^r$ by the complete
linear series of a line bundle $L$.  As in \cite{EKS}*{p.~514}, the line bundle
$L$ is called \define{determinantally presented} if the homogeneous ideal $I_{X
  \mid \PP^r}$ of $X$ in $\PP^r$ is generated by the $2$-minors of a
$1$\nobreakdash-generic matrix (i.e.{} no conjugate matrix has a zero entry) of
linear forms.  Definition~3.1 in \cite{G2} states that a property holds for a
\define{sufficiently ample} line bundle on $X$ if there exists a line bundle $A$
such that the property holds for all $L \in \Pic(X)$ for which $L \otimes
A^{-1}$ is ample.  Our main result is:

\begin{theorem}
  \label{t:main}
  Every sufficiently ample line bundle on a connected scheme is determinantally
  presented.
\end{theorem}

\noindent
We also describe, in terms of Castelnuovo-Mumford regularity, a collection of
determinantally presented line bundles on an arbitrary projective scheme; see
Corollary~\ref{c:eff}.

This theorem is a new incarnation of a well-known phenomenon --- roughly
speaking, the complexity of the first few syzygies of a projective subscheme is
inversely related to the positivity of the corresponding linear series.
Nevertheless, Theorem~\ref{t:main} counter-intuitively implies that most
projective embeddings by a complete linear series are simply the intersection of
a Segre variety with a linear subspace.  More precisely, if we fix the Euclidean
metric on the ample cone $\operatorname{Amp}(X)$ which it inherits from the
finite-dimensional real vector space $N^1(X) \otimes \RR$, then the fraction of
determinantally presented ample classes within distance $\rho$ of the trivial
class approaches $1$ has $\rho$ tends to $\infty$.

Theorem~\ref{t:main} also has consequences beyond showing that the homogeneous
ideal is generated by quadrics of rank at least $2$.  Proposition~6.13 in
\cite{Eisenbud} shows that an Eagon-Northcott complex is a direct summand of the
minimal graded free resolution of the ideal. Despite the classic examples, being
able to give a complete description of this resolution in the general setting
seems overly optimistic.  However, a determinantal presentation provides many
equations for higher secant varieties; see Proposition~1.3 in \cite{EKS}.  For a
scheme $X \subset \PP^r$, let $\operatorname{Sec}^k(X)$ be the Zariski closure
of the union of the linear spaces spanned by collections of $k+1$ points on $X$.
A natural generalization of Theorem~\ref{t:main} would be:

\begin{conjecture}
  \label{c:secants}
  Let $k$ be a positive integer.  If $X \subset \PP^r$ is embedded by the
  complete linear series of a sufficiently ample line bundle, then the
  homogeneous ideal of $\operatorname{Sec}^k(X)$ is generated by the
  $(k+2)$-minors of a $1$-generic matrix of linear forms.
\end{conjecture}

\noindent
This conjecture holds for rational normal curves (see Proposition~4.3 in
\cite{EisDet}), rational normal scrolls (see Proposition~2.2 in \cite{C-J}),
Segre varieties, and quadratic Veronese varieties (see \cite{SS}*{\S4}).  It
also extends the conjecture for curves appearing in \cite{EKS}*{p.~518} for
which \cite{Ravi} proves a set-theoretic version and \cite{Ginensky}*{\S7}
proves a scheme-theoretic version.  Although Theorem~1.1.4 in \cite{BGL}
produces counterexamples to this conjecture for some singular $X$,
Corollary~1.2.4 in \cite{BGL} provides supporting evidence when $X$ is smooth.
Theorem~1.3 in \cite{BB} suggests that the secant varieties in
Conjecture~\ref{c:secants} should be replaced by cactus varieties.

The secondary goal of this article is to effectively bound the determinantally
presented line bundles on specific schemes.  For an integral curve of genus $g$,
Theorem~1 in \cite{EKS} shows that a line bundle is determinantally presented
when its degree is at least $4g+2$ and this bound is sharp.  We provide the
analogous result on smooth varieties and Gorenstein toric varieties:

\begin{theorem}
  \label{t:effective}
  Let $X$ be a smooth variety of dimension $n$ or an $n$-dimensional Gorenstein toric variety
  and let $A$ be a very ample line bundle on $X$ such that $(X,A) \neq \bigl(
  \PP^n,\sO_\PP^n(1) \bigr)$.  If $B$ is a nef line bundle, $K_X$ is the
  dualizing bundle on $X$, and $L := K_X^2 \otimes A^j \otimes B$ with $j \geq
  2n+2$, then $L$ is determinantally presented.
\end{theorem}

\noindent
As an application of our methods, we describe determinantally presented ample
line bundles on products of projective spaces; see Theorem~\ref{t:proj}.

To prove these theorems, we need a source of appropriate matrices.  Composition
of linear series (a.k.a.{} multiplication in total coordinate ring or the Cox
ring) traditionally supply the required matrices. If $X \subset \PP^r$ is
embedded by the complete linear series for a line bundle $L$, then $H^0(X,L)$ is
the space of linear forms on $\PP^r$. Factoring $L$ as $L = E \otimes E'$ for
some $E, E' \in \Pic(X)$ yields a natural map $\mu \colon H^0(X,E) \otimes
H^0(X,E') \rTo H^0(X,E \otimes E') = H^0(X,L)$.  By choosing ordered bases $y_1,
\dotsc, y_s \in H^0(X,E)$ and $z_1, \dotsc, z_t \in H^0(X,E')$, we obtain an
associated $(s \times t)$\nobreakdash-matrix $\Omega := [ \mu(y_i \otimes z_j)
]$ of linear forms.  The matrix $\Omega$ is $1$-generic and its ideal
$\Ideal_2(\Omega)$ of $2$-minors vanishes on $X$; see Proposition~6.10 in
\cite{Eisenbud}.  Numerous classic examples of this construction can be found in
\cite{Room}.

With these preliminaries, the problem reduces to finding conditions on $E$ and
$E'$ which guarantee that $I_{X \mid \PP^r} = \Ideal_2(\Omega)$.  Inspired by
the approach in \cite{EKS}, Theorem~\ref{t:detPres} achieves this by placing
restrictions on certain modules arising from the line bundles $L$, $E$, and
$E'$. The key hypotheses require these modules to have a \define{linear free
  presentation}; the generators of the $\NN$-graded modules have degree $0$ and
their first syzygies must have degree $1$. Methods introduced by Green and
Lazarsfeld~\cites{G1,GL1} (for an expository account see \cite{Eisenbud}*{\S8},
\cite{G3}, or \cite{Lazarsfeld}*{\S1}) yield a cohomological criterion for our
modules to have a linear free presentation. Hence, we can prove
Theorem~\ref{t:main} by combining this with uniform vanishing results derived
from Castelnuovo-Mumford regularity. Building on known conditions (i.e.{}
sufficient conditions for a line bundle to satisfy $N_1$), we obtain effective
criteria for the appropriate modules to have a linear free presentation on
Gorenstein toric varieties, and smooth varieties.

Rather than focusing exclusively on a single factorization of the line bundle
$L$, we set up the apparatus to handle multiple factorizations; see
Lemma~\ref{l:diagram}.  Multiple factorizations of a line bundle were used in
\cite{GvBH} to study the equations and syzygies of elliptic normal curves and
their secant varieties.  They also provide a geometric interpretation for the
\emph{flattenings} appearing in \cite{GSS}*{\S7} and \cite{CGG}*{p.~1915}.
Using this more general setup, we are able to describe the homogeneous ideal for
every embedding of a product of projective spaces by a very ample line bundle as
the $2$-minors of appropriate $1$-generic matrices of linear forms; see
Proposition~\ref{p:mult}.

\subsection*{Conventions}
In this paper, $\NN$ is the set of nonnegative integers, $\id_W \in
\operatorname{Hom}(W,W)$ is the identity map, and $\bm{1} := (1, \dotsc, 1)$ is
the vector in which every entry is $1$.  We work over an algebraically closed
field $\kk$ of characteristic zero.  A variety is always
irreducible and all of our toric varieties are normal.  For a vector bundle $U$,
we write $U^j$ for the $j$-fold tensor product $U^{\otimes j} = U \otimes \dotsb
\otimes U$.

\subsection*{Acknowledgements}
We thank David Eisenbud, Tony Geramita, Rob Lazarsfeld, and Pete Vermeire for
helpful discussions.  The computer software \emph{Macaulay 2}~\cite{M2} was
useful for generating examples.  The first author was partially supported by NSF
grant DMS-0600471 and the Clare Boothe Luce program.  The second author was
partially supported by NSERC and grant KAW 2005.0098 from the Knut and Alice
Wallenberg Foundation.  We are grateful to the referee for the careful reading.

\section{Linear Free Presentations} 
\label{s:linsyz}

\noindent
This section collects the criteria needed to show that certain modules arising
from line bundles have a linear free presentation.  While accomplishing this, we
also establish some notation and nomenclature used throughout the document.

Let $X$ be a projective scheme over $\kk$, let $\sF$ be a coherent
$\sO_X$-module, and let $L$ be a line bundle on $X$.  We write $\Gamma(L) :=
H^0(X,L)$ for the $\kk$-vector space of global sections and $S := \Sym \bigl(
\Gamma(L) \bigr)$ for the homogeneous coordinate ring of $\PP^r := \PP\bigl(
\Gamma(L) \bigr)$.  Consider the $\NN$\nobreakdash-graded $S$\nobreakdash-module
$F := \bigoplus_{j \geq 0} H^0(X, \sF \otimes L^j)$.  When $\sF = \sO_X$, $F$ is
the section ring of $L$.  However, when $\sF = L$, the module $F$ is the
truncation of the section ring omitting the zeroth graded piece and shifting
degrees by $-1$.  Let $P_\bullet$ be a minimal graded free resolution of $F$:
\[
\begin{diagram}[h=1.3em,w=1.5em]
  \dotsb & \rTo & & \bigoplus S(-a_{i,j}) & & \rTo & \dotsb & \rTo & & \bigoplus
  S(-a_{1,j}) & & \rTo & & \bigoplus S(-a_{0,j}) & & \rTo & F & \rTo & 0 \, . \\
  & & & \dEqual & & & & & & \dEqual & & & & \dEqual & & & & \\
  & & & P_i & & & & & & P_1 & & & & P_0
\end{diagram} 
\]
Following \cite{EKS}*{p.~515}, we say that, for $p \in \NN$, $\sF$ has a
\define{linear free resolution to stage $p$ with respect to $L$} or $F$ has a
\define{linear free resolution to stage~$p$} if $P_i = \bigoplus S(-i)$ for all
$0 \leq i \leq p$.  Thus, $F$ has a linear free resolution to stage $0$ if and
only if it is generated in degree $0$.  Since having a linear free resolution to
stage~$1$ implies that the relations among the generators (a.k.a.{} first
syzygies) are linear, the module $F$ has a linear free resolution to stage~$1$
if and only if it has a linear free presentation.  In this case, we say that
$\sF$ has a \define{linear free presentation with respect to $L$}.  More
generally, having a linear free resolution to stage~$p$ is the module-theoretic
analogue of the $N_p$-property introduced in \cite{GL1}*{\S3}.  If $X$ is
connected, then the line bundle $L$ satisfies $N_1$ precisely when $L$ has a
linear free presentation with respect to itself and satisfies $N_p$ when $L$ has
a linear free resolution to stage $p$.  Following Convention~0.4 in \cite{EL},
we do not assume that $X$ is normal.

Henceforth, we assume that $L$ is globally generated.  In other words, the
natural evaluation map $\operatorname{ev}_{L} \colon \Gamma(L) \otimes_\kk \sO_X
\rTo L$ is surjective. If $M_{L} := \Ker(\operatorname{ev}_{L})$, then $M_{L}$
is a vector bundle of rank $r := \dim_\kk \Gamma(L) -1$ which sits in the short
exact sequence
\begin{equation}
  \tag{$\ast$}
  \label{e:ses}
  \begin{diagram}[w=1.3em]
    0 & & \rTo & M_L & & \rTo & & \Gamma(L) \otimes_\kk \sO_X  & &
    \rTo & & L & \rTo & & 0 \, .
  \end{diagram}
  \addtocounter{equation}{1}
\end{equation}
For convenience, we record the following cohomological criteria which is a minor
variant of Theorem~5.6 in \cite{Eisenbud}, Proposition~2.4 in \cite{G3}, or
Lemma~1.6 in \cite{EL}.

\begin{lemma}
  \label{l:criteria}
  If $H^1(X, \textstyle\bigwedge^i M_{L} \otimes \sF \otimes L^j) = 0$ for all
  $1 \leq i \leq p+1$ and all $j \geq 0$, then the coherent $\sO_X$-module $\sF$
  has a linear free resolution to stage $p$ with respect to $L$.  In
  characteristic zero, $\bigwedge^i M_{L}$ is a direct summand of $M_{L}^i$, so
  it suffices to show $H^1(X, M_{L}^i \otimes \sF \otimes L^j) = 0$ for all $1
  \leq i \leq p+1$ and all $j \geq 0$.  \qed
\end{lemma}

\begin{proof}[Sketch of Proof]
  The key observation is that the graded Betti numbers for the minimal free
  resolution of $F$ can be computed via Koszul cohomology.  If $L$ is globally
  generated and $\PP^r = \PP \bigl( H^0(X, L) \bigr)$, then there is a morphism
  $\varphi_L \colon X \rTo \PP^r$ with $\varphi_L^* \bigl(\sO_{\PP^r}(1) \bigr) =
  L$.  Since the pullback by $\varphi_L^*$ of $0 \rTo M_{\sO_{\PP^r}(1)} \rTo
  \Gamma\bigl( \sO_{\PP^r}(1) \bigr) \otimes_\kk \sO_{\PP^r} \rTo \sO_{\PP^r}(1)
  \rTo 0$ is just \eqref{e:ses}, the proof of Theorem~5.6 in \cite{Eisenbud}
  goes through working on $X$ instead of $\PP^r$.
\end{proof}

Multigraded Castelnuovo-Mumford regularity, as developed in \cite{MS}*{\S6} or
\cite{HSS}*{\S2}, allows us to exploit this criteria.  To be more precise, fix a
list $B_1, \dotsc, B_\ell$ of globally generated line bundles on $X$.  For a
vector $\bm{u} := (u_1, \dotsc, u_\ell) \in \ZZ^\ell$, we set $\bm{B}^{\bm{u}}
:= B_1^{u_1} \otimes \dotsb \otimes B_\ell^{u_\ell}$ and we write $\fB := \{
\bm{B}^{\bm{u}} : \bm{u} \in \NN^{\ell} \} \subset \operatorname{Pic}(X)$ for
the submonoid generated by these line bundles.  If $\bm{e}_1, \dotsc,
\bm{e}_\ell$ is the standard basis for $\ZZ^\ell$ then $\bm{B}^{\bm{e}_j} =
B_j$.  A coherent $\sO_X$-module $\sF$ is said to be \define{regular with
  respect to $B_1, \dotsc, B_\ell$} if $H^i(X, \sF \otimes \bm{B}^{-\bm{u}}) =
0$ for all $i > 0$ and all $\bm{u} \in \NN^\ell$ satisfying $|\bm{u}| := u_1 +
\dotsb + u_\ell = i$.  When $\ell = 1$, we recover the version of
Castelnuovo-Mumford regularity found in \cite{PAG}*{\S1.8}.

Although the definition may not be intuitive, the next result shows that regular
line bundles are at least ubiquitous.

\begin{lemma}
  \label{l:lower}
  Let $X$ be a scheme and let $B_1, \dotsc, B_\ell$ be globally generated line
  bundles on $X$.  If there is a positive vector $\bm{w} \in \ZZ^\ell$ such that
  $\bm{B}^{\bm{w}}$ is ample, then a sufficiently ample line bundle on $X$ is
  regular with respect to $B_1, \dotsc, B_\ell$.
\end{lemma}

\noindent
The hypothesis on $\bm{w}$ means that the cone $\pos(B_1, \dotsc, B_\ell)$
generated by $B_1, \dotsc, B_\ell$ contains an ample line bundle.  In other
words, the subcone $\pos(B_1, \dotsc, B_\ell)$ of $\Nef(X)$ has a nonempty
intersection with the interior of $\Nef(X)$.

\begin{proof}
  It suffices to find a line bundle $A$ on $X$ such that, for any nef line
  bundle $C$, $A \otimes C$ is regular with respect to $B_1, \ldots, B_{\ell}$.
  Because $\bm{B}^{\bm{w}}$ is ample, Fujita's Vanishing Theorem (e.g.{}
  Theorem~1 in \cite{Fujita}) implies that there is $k \in \NN$ such that, for any
  nef line bundle $C$, we have $H^i(X, \bm{B}^{j \bm{w}} \otimes C) = 0$ for all
  $i > 0$ and all $j \geq k$.  Let $n := \dim X$ and consider $A :=
  \bm{B}^{(k+n)\bm{w}}$.  Since $\bm{w}$ is positive, the line bundle
  $\bm{B}^{n\bm{w} - \bm{u}}$ is nef for all $\bm{u} \in \NN^{\ell}$ with $0
  \leq |\bm{u}| \leq n$.  Therefore, we have $H^i(X, (A \otimes C) \otimes
  \bm{B}^{-\bm{u}}) = H^i (X, \bm{B}^{k \bm{w}} \otimes (\bm{B}^{n
    \bm{w}-\bm{u}} \otimes C)) = 0$ for all $i > 0$ and all $\bm{u} \in
  \NN^{\ell}$ satisfying $|\bm{u}| = i$.
\end{proof}

Before describing the pivotal results in this section, we record a technical
lemma bounding the regularity of certain tensor products.  Our approach is a
hybrid of Proposition~1.8.9 and Remark~1.8.16 in \cite{PAG}.

\begin{lemma}
  \label{l:regPres}
  Let $X$ be a scheme of dimension $n$ and $\sF$ be a coherent $\sO_X$-module.
  Fix a vector bundle $V$ and a globally generated ample line bundle $B$ on $X$.
  If $m$ is positive integer such that $\sF$, $V$, and $B^m$ are all regular
  with respect to $B$, then $\sF \otimes V \otimes B^w$ is also regular with
  respect to $B$ for all $w \geq (m-1)(n-1)$.
\end{lemma}

\begin{proof}
  Since $\sF$ and $B^m$ are regular with respect to $B$, Corollary~3.2 in
  \cite{arapura} or Theorem~7.8 in \cite{MS} (cf.{} Proposition~1.8.8 in
  \cite{PAG}) produces a locally free resolution of $\sF$ of the form 
  $
  \begin{diagram}[w=1.0em]
    \dotsb & & \rTo & & \bigoplus B^{-jm} & & \rTo & & \dotsb & &
    \rTo & & \bigoplus B^{-m} & & \rTo & & \bigoplus \sO_X & &
    \rTo & & \sF & & \rTo & 0 \, .
  \end{diagram}
  $ 
  Tensoring by a locally free sheaf perserves exactness, so we obtain the
  exact complex
  \[
  \begin{diagram}[w=1.3em]
    \dotsb & & \rTo & & \bigoplus V \otimes B^{w-jm} & & \rTo & & \dotsb & &
    \rTo & & \bigoplus V \otimes B^{w} & & \rTo & & \sF \otimes V \otimes B^{w}
    & & \rTo & 0 \, .
  \end{diagram}
  \]
  Since $V$ is also regular with respect to $B$, Mumford's Lemma (e.g.{}
  Theorem~1.8.5 in \cite{PAG}) implies that $H^{i+j}(X, V \otimes B^{w-jm-i}) =
  0$ for $i \geq 1$ provided we have $w -jm-i \geq -i-j$.  Chasing through the
  complex (see Proposition~B.1.2 in \cite{PAG}), we conclude that $\sF \otimes V
  \otimes B^w$ is also regular with respect to $B$ when $w \geq (m-1)(n-1)$.
\end{proof}

The following three propositions each provide sufficient conditions for an
appropriate line bundle to have a linear free presentation with respect to
another line bundle.

\begin{proposition}
  \label{p:lfp}
  Fix a positive integer $m$ and a scheme $X$ of dimension $n$.  Let $L$ be a
  line bundle on $X$ and let $B$ be a globally generated ample line bundle on
  $X$.  If $L^j$ and $B^m$ are regular with respect to $B$ for all $j \geq 1$,
  then $B^w$ has a linear free presentation with respect to $L$ for all $w \geq
  2(m-1)n+1$.
\end{proposition}

\begin{proof}
  We first prove that $M_{L} \otimes B^m$ is regular with respect to $B$.
  Tensoring \eqref{e:ses} with $B^{m - i}$ and taking the associated long exact
  sequence gives
  \begin{multline*}
    \begin{diagram}[w=1.3em]
      \Gamma(L) \otimes H^0(X, B^{m - i}) & & \rTo & & H^0(X, L \otimes B^{m -
        i}) & & \rTo & & H^1(X,M_{L} \otimes B^{m - i}) & & \rTo & \dotsb
    \end{diagram} \\
    \begin{diagram}[w=1.3em]
      & \rTo & & H^{i-1}(X, L \otimes B^{m - i}) & & \rTo & & H^i(X,M_{L}
      \otimes B^{m - i}) & & \rTo & & \Gamma(L) \otimes H^i(X, B^{m -i}) \, .
    \end{diagram}
  \end{multline*}
  Since $L$ is regular with respect to $B$, Mumford's Lemma (e.g.{}
  Theorem~1.8.5 in \cite{PAG}) shows that, for all $k \in \NN$, the map
  $\Gamma(L) \otimes H^{0}(X,B^k) \rTo H^0(X, L \otimes B^k)$ is surjective and,
  for all $i > 0$ and all $k \in \NN$, we have $H^i(X, L \otimes B^{k-i}) = 0$.
  As $m$ is a positive integer, the map $\Gamma(L) \otimes H^{0}(X,B^{m-1}) \rTo
  H^0(X, L \otimes B^{m-1})$ is surjective and $H^{i-1}(X, L \otimes B^{m-i}) =
  0$ for all $i > 1$.  Since $B^m$ is also regular with respect to $B$, we have
  $H^i(X, B^{m-i}) = 0$ for all $i > 0$.  It follows that $H^i(X,M_{L} \otimes
  B^{m-i}) = 0$ for all $i > 0$.

  By Lemma~\ref{l:criteria}, it suffices to show that, for all $j \in \NN$, we
  have $H^1(X, M_{L} \otimes B^{w}\otimes L^j) = 0$ and $H^1(X, M_{L}^2 \otimes
  B^{w} \otimes L^j) = 0$.  Thus, it suffices to show that the vector bundles
  $M_{L} \otimes B^{w+1} \otimes L^j$ and $M_{L}^{2} \otimes B^{w+1} \otimes
  L^j$ are both regular with respect to $B$.  If $w \geq (m-1)n$, then
  Lemma~\ref{l:regPres} implies that $(M_L \otimes B^m) \otimes L^j \otimes
  B^{w+1-m} = M_L \otimes B^{w+1} \otimes L^j$ is regular with respect to $B$.
  Similarly, if $w \geq 2(m-1)n+1$, then using Lemma~\ref{l:regPres} twice
  establishes that the vector bundle
  \[
  \bigl( (M_L \otimes B^m) \otimes (M_L \otimes B^m) \otimes B^{(m-1)(n-1)}
  \bigr) \otimes L^j \otimes B^{w-mn-m+n} = M_{L}^2 \otimes B^{w+1} \otimes L^j
  \]
  is also regular with respect to $B$.
\end{proof}

By adapting the proof of Theorem~1.1 in \cite{HSS}, we obtain the second
proposition.

\begin{proposition}
  \label{p:regCrit}
  Let $\bm{m} \in \NN^{\ell}$ be a vector satisfying $\bm{B}^{\bm{m}- \bm{e}_j}
  \in \fB$ for all $1 \leq j \leq \ell$ and let the coherent $\sO_X$-module
  $\sF$ be regular with respect to $B_1, \dotsc, B_\ell$.  If $L :=
  \bm{B}^{\bm{m}}$ and the map
  \[
  \begin{diagram}[w=1.3em]
    \Gamma(L)\otimes H^0(X,\sF \otimes \bm{B}^{-\bm{e}_j}) & & \rTo & & H^{0}(X,
    \sF \otimes \bm{B}^{\bm{m}-\bm{e}_j})
  \end{diagram}
  \] 
  is surjective for all $1 \leq j \leq \ell$, then $\sF$ has a linear
  presentation with respect to $L$.
\end{proposition}

\noindent
The condition that $\bm{B}^{\bm{m} - \bm{e}_j} \in \fB$ for all $1 \leq j \leq
\ell$ implies that $L = \bm{B}^{\bm{m}}$ lies in the interior of the cone
$\pos(B_1, \dotsc, B_\ell)$.

\begin{proof}
  We first prove that $M_{L} \otimes \sF$ is regular with respect to $B_1,
  \dotsc, B_\ell$.  Tensoring \eqref{e:ses} with $\sF \otimes \bm{B}^{-\bm{u}}$
  and taking the associated long exact sequence gives
  \begin{multline*}
    \begin{diagram}[w=1.2em]
      \Gamma(L) \otimes H^0(X, \sF \otimes \bm{B}^{-\bm{u}}) & & \rTo & & H^0(X,
      \sF \otimes \bm{B}^{\bm{m} - \bm{u}}) & & \rTo & & H^1(X,M_{L} \otimes \sF
      \otimes \bm{B}^{-\bm{u}}) & & \rTo & \dotsb
    \end{diagram} \\
    \begin{diagram}[w=1.2em]
      & \rTo & & H^{i-1}(X, \sF \otimes \bm{B}^{\bm{m} - \bm{u}}) & & \rTo & &
      H^i(X,M_{L} \otimes \sF \otimes \bm{B}^{-\bm{u}}) & & \rTo & & \Gamma(L)
      \otimes H^i(X, \sF \otimes \bm{B}^{-\bm{u}}) \, .
    \end{diagram}
  \end{multline*}
  Since $\sF$ is regular with respect to $B_1, \dotsc, B_\ell$, Theorem~2.1 in
  \cite{HSS} shows that, for all $i > 0$ and all $\bm{u}, \bm{v} \in \NN^{\ell}$
  with $|\bm{u}| = i$, $H^i(X, \sF \otimes \bm{B}^{\bm{v}-\bm{u}}) = 0$.  As
  $\bm{B}^{\bm{m}-\bm{e}_j} \in \fB$ for $1 \leq j \leq \ell$, we see that
  $H^{i-1}(X, \sF \otimes \bm{B}^{\bm{m}-\bm{u}}) = 0$ for all $i > 1$ and all
  $\bm{u} \in \NN^{\ell}$ satisfying $|\bm{u}| = i$.  By hypothesis, the map
  $\Gamma(L) \otimes H^0(X, \sF \otimes \bm{B}^{-\bm{e}_j}) \rTo H^{0}(X, \sF
  \otimes \bm{B}^{\bm{m}-\bm{e}_j})$ is surjective for all $1 \leq j \leq \ell$.
  It follows that $H^i(X,M_{L} \otimes \sF \otimes \bm{B}^{-\bm{u}}) = 0$ for
  all $i > 0$ and all $\bm{u} \in \NN^\ell$ such that $|\bm{u}| = i$.

  By Lemma~\ref{l:criteria}, it suffices to show that $H^1(X, M_{L} \otimes \sF
  \otimes L^j)$ and $H^1(X, M_{L}^2 \otimes \sF \otimes L^j)$ are zero for $j
  \in \NN$.  Since $M_{L} \otimes \sF$ is regular with respect to $B_1, \dotsc,
  B_\ell$, the vanishing of the first group follows from Theorem~2.1 in
  \cite{HSS}.  For the second, tensoring \eqref{e:ses} with $M_{L} \otimes \sF
  \otimes L^j$ gives the exact sequence:
  \[
  \begin{diagram}[w=1.1em]
    \Gamma(L) \otimes H^0(X,M_{L} \otimes \sF \otimes L^j) & & \rTo & & H^0(X,
    M_{L} \otimes \sF \otimes L^{j+1}) & & \rTo & & H^1(X,M_{L}^2 \otimes \sF
    \otimes L^j) & & \rTo & 0 \, .
  \end{diagram}
  \]
  Because $M_{L} \otimes \sF$ is regular with respect to $B_1, \dotsc, B_\ell$,
  Theorem~2.1 in \cite{HSS} also shows that the left map is surjective for all
  $j \geq 0$.
\end{proof}

Our third proposition is a variant of Proposition~3.1 in \cite{EL}.

\begin{proposition}
  \label{p:smoothAdj}
  Let $X$ be a smooth variety of dimension $n$, let $K_X$ be its canonical
  bundle, and let $A$ be a very ample line bundle on $X$ such that $(X,A) \neq
  \bigl( \PP^n, \sO_{\PP^n}(1) \bigr)$.  Suppose that $B$ and $C$ are nef line
  bundles on $X$. If the integers $w$ and $m$ are both greater than $n$, then
  the line bundle $K_X \otimes A^w \otimes B$ has a linear free presentation
  with respect to $K_X \otimes A^m \otimes C$.
\end{proposition}

\begin{proof}
  Let $\sF := K_X \otimes A^w \otimes B$ and $L := K_X \otimes A^m \otimes C$.
  Since Proposition~3.1 in \cite{EL} shows that $L$ satisfies $N_0$ and
  Equation~3.2 in \cite{EL} shows that $H^1(X, M_{L}^i \otimes \sF \otimes L^j)
  = 0$ for all $1 \leq i \leq 2$ and all $j \geq 0$, Lemma~\ref{l:criteria}
  completes the proof.
\end{proof}

\section{Determinantally Presented Line Bundles}
\label{s:minors}

\noindent
The goal of this section is to prove Theorem~\ref{t:main}.  We realize this goal
by developing general methods for showing that a line bundle is determinantally
presented; see Theorem~\ref{t:detPres}.

Suppose that $X \subset \PP^r$ is embedded by the complete linear series for a
line bundle $L$.  Factor $L$ as $L = E \otimes E'$ for some $E, E' \in \Pic(X)$
and let $\mu_{E,E'} \colon H^0(X,E) \otimes H^0(X,E') \rTo H^0(X,L)$ denote the
natural multiplication map. Choose ordered bases $y_1, \dotsc, y_s$ and $z_1,
\dotsc, z_t$ for the $\kk$\nobreakdash-vector spaces $H^0(X,E)$ and $H^0(X,E')$
respectively. Define $\Omega = \Omega(E,E')$ to be the associated $(s \times
t)$-matrix $[\mu_{E,E'}(y_i \otimes z_j)]$ of linear forms.  Its ideal
$\Ideal_2(\Omega)$ of $2$-minors is independent of the choice of bases.
Proposition~6.10 in \cite{Eisenbud} shows that $\Omega$ is $1$-generic and that
$\Ideal_2(\Omega)$ vanishes on $X$.

Inspired by \cite{EKS}*{\S2}, the key technical result is:

\begin{lemma} 
  \label{l:diagram}
  If $L$ is a very ample line bundle on $X$ satisfying $N_1$ and $\{ (E_i^{\,},
  E'_i) \}$ is a family of factorizations for $L$, then the commutative diagram
  \eqref{e:diagram} has exact rows and columns.  Moreover, if $\varphi_{\, 2}$
  is surjective, then the homogeneous ideal $I_{X \mid \PP^r}$ is generated by
  the $2$-minors of the matrices $\Omega(E_i^{\,},E'_i)$ if and only if $Q_2$
  surjects onto $Q_1$.
\end{lemma}

\begin{figure}[ht]
  \begin{equation}
    \tag{$\maltese$}
    \label{e:diagram}
    \begin{diagram}[h=1em,w=2.35em]
      & & & & 0 & & & & 0 & &\\
      & & & & \dTo & & & & \dTo \\
      & & & & \bigoplus\nolimits_i \textstyle\bigwedge^{2}
      \Gamma(E_i^{\,}) \otimes \textstyle\bigwedge^{2}
      \Gamma(E'_i) & & \rTo^{\varphi} & & (I_{X \mid \PP^r})^{\,}_2 & & \\
      & & & & & & & & & & \\
      & & & & \dTo & & & & \dTo & &\\
      & & & & & & & & & & \\
      0 & \rTo & Q_2 & \rTo & \bigoplus\nolimits_i \Sym_2 \bigl(
      \Gamma(E_i^{\,}) \otimes \Gamma(E'_i) \bigr) & &
      \rTo^{\varphi_{\, 2}} & & \Sym_2 \bigl( \Gamma(L) \bigr) & & \\
      & & & & & & & & & &\\
      & & \dTo^{\psi} & & \dTo & & & & \dTo & &\\
      & & & & & & & & & & \\
      0 & \rTo & Q_1 & \rTo & \bigoplus\nolimits_i \Sym_2 \bigl(
      \Gamma(E_i^{\,}) \bigr) \otimes \Sym_2 \bigl( \Gamma(E'_i)
      \bigr) & & \rTo^{\varphi_1} & & \Gamma(L^2) & & \\
      & & & & & & & & & &\\
      & & & & \dTo & & & & \dTo & &\\
      & & & & 0 & & & & 0 & & \;
    \end{diagram}
  \end{equation}
  \addtocounter{equation}{1}
\end{figure}

\begin{proof}
  To begin, we prove the columns are exact.  Since $L$ satisfies $N_0$ (i.e.{}
  the natural maps $\Sym_j \bigl(\Gamma(L) \bigr) \rTo H^0(X,L^j)$ are
  surjective for all $j \in \NN$), the ideal $I_{X \mid \PP^r}$ is the kernel of
  the map from the homogeneous coordinate ring of $\PP^r$ to the section ring of
  $L$.  By taking the quadratic components, we obtain the right column.  The
  middle column is the direct sum of the complexes:
  \[
  \begin{diagram}[w=1em]
    0 & \rTo & \textstyle\bigwedge^2 \Gamma(E_i^{\,}) \otimes
    \textstyle\bigwedge^2 \Gamma(E'_i) & \rTo & \Sym_2 \bigl( \Gamma(E_i^{\,})
    \otimes \Gamma(E'_i) \bigr) & \rTo & \Sym_2 \bigl( \Gamma(E_i^{\,}) \bigr)
    \otimes \Sym_2 \bigl( \Gamma( E'_i) \bigr) & \rTo & 0 .
  \end{diagram}
  \]
  The map $\textstyle\bigwedge^2 \Gamma(E_i^{\,}) \otimes \textstyle\bigwedge^2
  \Gamma(E'_i) \rTo \Sym_2 \bigl( \Gamma(E_i^{\,}) \otimes \Gamma(E'_i) \bigr)$,
  defined by
  \begin{equation}
    \tag{$\dagger$}
    \label{e:det}
    \begin{diagram}[w=1.5em]
      e \wedge f \otimes e' \wedge f' & & \rMapsto & & (e \otimes e') \cdot
      (f\otimes f') - (e\otimes f') \cdot (f \otimes e') \, ,
    \end{diagram}
    \addtocounter{equation}{1}
  \end{equation} 
  is simply the inclusion map determined by the $2$-minors of the generic
  matrix.  The map $\Sym_2 \bigl( \Gamma(E_i^{\,}) \otimes \Gamma(E'_i) \bigr)
  \rTo \Sym_2 \bigl( \Gamma(E_i^{\,}) \bigr) \otimes \Sym_2 \bigl( \Gamma( E'_i)
  \bigr)$ is $(e\otimes e') \cdot (f \otimes f') \rMapsto ef \otimes e'f'$.
  Hence, each of these complexes is exact, so the middle column also is.  By
  definition, $Q_1$ and $Q_2$ are the kernels of the $\varphi_1$ and
  $\varphi_{\, 2}$ respectively, and $\psi$ is the induced map between them.

  We next identify the horizontal maps.  By applying the functor $\Sym_2$ to
  $\mu_{E_i^{\,}, E'_i}$, we obtain a map from $\Sym_2 \bigl( \Gamma(E_i^{\,})
  \otimes \Gamma(E'_i) \bigr)$ to $\Sym_2 \bigl( \Gamma(L) \bigr)$ for each $i$;
  $\varphi_{\, 2}$ is their direct sum.  The composite map $\mu_{L,L} \circ (
  \mu_{E_i^{\,},E'_i} \otimes \mu_{E_i^{\,},E'_i}) \colon \Gamma(E_i^{\,})
  \otimes \Gamma(E'_i) \otimes \Gamma(E_i^{\,}) \otimes \Gamma(E'_i) \rTo
  \Gamma(L^2)$ factors through $\Sym_2 \bigl( \Gamma(E_i^{\,}) \bigr) \otimes
  \Sym_2 \bigl( \Gamma(E'_i) \bigr)$, and $\varphi_1$ is the direct sum of the
  associated maps from $\Sym_2\bigl( \Gamma(E_i^{\,}) \bigr) \otimes \Sym_2
  \bigl( \Gamma(E'_i) \bigr)$ to $\Gamma(L^2)$.  The map $\varphi$ is induced by
  $\varphi_{\, 2}$.  From \eqref{e:det}, we see that the image of $\varphi$ is
  generated by the $2$-minors of the matrices $\Omega(E_i^{\,},E'_i)$.

  Finally, the line bundle $L$ satisfies $N_1$, so the quadratic component
  $(I_{X \mid \PP^r})_2$ generates the entire ideal $I_{X \mid \PP^r}$.  Hence,
  the image of $\varphi$ generates the ideal $I_{X \mid \PP^r}$ if and only if
  $\varphi$ is surjective.  Since $\varphi_2$ is surjective, the Snake Lemma
  (e.g.{} Lemma~1.3.2 in \cite{Weibel}) shows that the surjectivity of $\varphi$
  is equivalent to the surjectivity of $\psi$.
\end{proof}

Our main application for Lemma~\ref{l:diagram} focuses on a single factorization
of the line bundle $L$.  The proof follows the strategy in \cite{EKS}*{\S2}.

\begin{theorem}
  \label{t:detPres}
  Let $L$ be a very ample line bundle on a scheme $X$ satisfying $N_1$.  If $L =
  E \otimes E'$ for some nontrivial $E, E' \in \Pic(X)$ and the following
  conditions hold
  \begin{enumerate}
  \item $E$ has a linear presentation with respect to $E'$,
  \item $E'$ has a linear presentation with respect to $E$,
  \item $E^2$ has a linear presentation with respect to $E'$,
  \item both $E$ and $E'$ satisfy $N_0$,
  \end{enumerate}
  then the $2$-minors of the matrix $\Omega(E,E')$ generate the homogeneous
  ideal of $X$ in $\PP\bigl( \Gamma(L) \bigr)$.  In particular, the line bundle
  $L$ is determinantally presented.
\end{theorem}
 
\begin{proof}
  Given Lemma~\ref{l:diagram}, it suffices to show that the map $\psi \colon Q_2
  \rTo Q_1$ is surjective.  To accomplish this, we reinterpret both modules.
  Since Condition~(a) or (b) imply that the map $\mu_{E,E'} \colon \Gamma(E)
  \otimes \Gamma(E') \rTo \Gamma(L)$ is surjective, we obtain an exact sequence
  \[
  \begin{diagram}[w=1.5em]
    \Ker(\mu_{E,E'}) \otimes \Gamma(E) \otimes \Gamma(E') & & \rTo & & \Sym_2
    \bigl( \Gamma(E) \otimes \Gamma(E') \bigr) & & \rTo & & \Sym_2 \bigl(
    \Gamma(L) \bigr) & \rTo & & 0 \, ,
  \end{diagram}
  \]
  so the image of $\Ker(\mu_{E,E'}) \otimes \Gamma(E) \otimes \Gamma(E')$
  generates $Q_2$ in $\Sym_2 \bigl( \Gamma(E) \otimes \Gamma(E') \bigr)$.  The
  maps $\mu_{E,E}$ and $\mu_{E',E'}$ factor through $\Sym_2 \bigl( \Gamma(E)
  \bigr)$ and $\Sym_2 \bigl( \Gamma(E') \bigr)$ and thus induce maps $\eta
  \colon \Sym_2 \bigl( \Gamma(E) \bigr) \rTo \Gamma(E^2)$ and $\eta' \colon
  \Sym_2 \bigl( \Gamma(E') \bigr) \rTo \Gamma({E'}^2)$ respectively.  It follows
  that $\varphi_1$ is the composition $\mu_{E^2,{E'}^2} \circ (\eta \otimes
  \eta') \colon \Sym_2 \bigl( \Gamma(E) \bigr) \otimes \Sym_2 \bigl( \Gamma(E')
  \bigr) \rTo \Gamma(E^2 \otimes E'^2) = \Gamma(L^2)$.  Hence, $Q_1$ is the sum
  of the images of $\Ker(\eta) \otimes \Gamma(E') \otimes \Gamma(E')$ and
  $\Gamma(E) \otimes \Gamma(E) \otimes \Ker(\eta')$, and the pullback to $\Sym_2
  \bigl( \Gamma(E) \bigr) \otimes \Sym_2 \bigl( \Gamma(E') \bigr)$ of
  $\Ker(\mu_{E^2,E'^2})$.

  We now break the proof that $Q_2$ surjects onto $Q_1$ into four steps:
  \begin{enumerate}
  \item[(i)] the image of $\Ker(\mu_{E^2, E'}) \otimes \Gamma(E')$ in
    $\Gamma(E^2) \otimes \Gamma(E'^2)$ contains $\Ker(\mu_{E^2, E'^2})$,
  \item[(ii)] the image of $\Ker(\mu_{E, E'}) \otimes \Gamma(E)$ in $\Gamma(E^2)
    \otimes \Gamma(E')$ contains $\Ker(\mu_{E^2, E'})$,
  \item[(iii)] the image of $\Ker(\mu_{E, E'}) \otimes \Gamma(E)$ in $\Sym_2
    \bigl( \Gamma(E) \bigr) \otimes \Gamma(E')$ contains $\Ker(\eta) \otimes
    \Gamma(E')$,
  \item[(iv)] the image of $\Ker(\mu_{E, E'}) \otimes \Gamma(E')$ in $\Gamma(E)
    \otimes \Sym_2 \bigl( \Gamma(E') \bigr)$ contains $\Gamma(E) \otimes
    \Ker(\eta')$.
  \end{enumerate}
  Tensoring with $\kk$-vector space $\Gamma(E')$, Step~(ii) yields a surjective
  map
  \[
  \begin{diagram}[w=1.5em]
    \Ker(\mu_{E, E'}) \otimes \Gamma(E) \otimes \Gamma(E') & & \rTo & &
    \Ker(\mu_{E^2, E'}) \otimes \Gamma(E') \, .
  \end{diagram}
  \] 
  Combining this with Step~(i) shows that the map $\Ker(\mu_{E, E'}) \otimes
  \Gamma(E) \otimes \Gamma(E') \rTo \Ker(\mu_{E^2, E'^2})$ is surjective.  Again
  by tensoring with $\kk$-vector space $\Gamma(E')$, Step~(iii) gives a
  surjective map $\Ker(\mu_{E, E'}) \otimes \Gamma(E) \otimes \Gamma(E') \rTo
  \Ker(\eta) \otimes \Gamma(E') \otimes \Gamma(E')$.  Similarly, Step~(iv)
  implies that the map $\Ker(\mu_{E, E'}) \otimes \Gamma(E) \otimes \Gamma(E')
  \rTo \Gamma(E) \otimes \Gamma(E) \otimes \Ker(\eta')$ is surjective.
  Therefore, it is enough to establish the four steps.

  For Step~(i), Condition~(c) implies that $\Ker(\mu_{E^2, E'})$, the span of
  the linear relations on $\bigoplus_{j \geq 0} H^0(X, E^2 \otimes E'^j)$
  regarded as a $\Sym \bigl( \Gamma(E') \bigr)$-module, generates the relations
  in higher degrees as well.  Hence, $\Ker(\mu_{E^2, E'}) \otimes \Gamma(E')$
  maps onto the quadratic relations which are the kernel of the composite map
  $\mu_{E^2,E'^2} \circ ( \id_{\Gamma(E^2)} \otimes \eta')$.  Since this kernel
  is generated by $\Gamma(E^2) \otimes \Ker(\eta')$ and the pullback of
  $\Ker(\mu_{E^2, E^{'2}})$, Condition~(d) implies that $\eta'$ is surjective,
  and we have established Step~(i).

  To complete the proof, we simultaneously establish Steps~(ii) and (iii); the
  symmetric argument yields Step~(iv).  Condition~(b) implies that $\Ker(\mu_{E,
    E'})$ generates all the relations on $\bigoplus_{j \geq 0} H^0(X, E' \otimes
  E^j)$ regarded as a $\Sym \bigl( \Gamma(E) \bigr)$-module.  In particular, the
  vector space $\Ker(\mu_{E, E'}) \otimes \Gamma(E)$ maps onto the quadratic
  relations which are the kernel of the composite map $\mu_{E^2,E'} \circ ( \eta
  \otimes \id_{\Gamma(E')})$.  This kernel is generated by $\Ker(\eta) \otimes
  \Gamma(E')$ and the pullback of $\Ker(\mu_{E^2, E'})$.  Condition~(d) implies
  that $\eta$ is surjective, so Step~(ii) and Step~(iii) follow.
\end{proof}

As the proof indicates, Theorem~\ref{t:detPres} holds under a weaker version of
Condition~(d).  Specifically, it is only necessary that $\eta$ and $\eta'$ are
surjective.  Nevertheless, in all of our applications, a stronger condition is
satisfied: both $E$ and $E$ satisfy $N_1$.
 
This theorem leads to a description, given in terms of Castelnuovo-Mumford
regularity, for some determinantally presented line bundles on any projective
scheme.

\begin{corollary}
  \label{c:eff}
  Let $X$ be a connected scheme and let $B_1, \dotsc, B_\ell$ be globally
  generated line bundles on $X$ for which there exists $\bm{w} \in \NN^\ell$
  such that $\bm{B}^{\bm{w}}$ is ample.  If $\bm{B}^{\bm{m}}$ is regular with
  respect to $B_1, \dotsc, B_\ell$ for $\bm{m} \in \NN^\ell$ and
  $\bm{B}^{2\bm{m}}$ is very ample, then the line bundle
  $\bm{B}^{2\bm{m}+\bm{u}}$ is determinantally presented for any $\bm{u} \in
  \NN^{\ell}$.
\end{corollary}

\begin{proof}
  Factor $L := \bm{B}^{2\bm{m}+\bm{u}}$ as $L = E \otimes E'$ where $E :=
  \bm{B}^{\bm{m}}$ and $E' := \bm{B}^{\bm{m}+\bm{u}}$.  Theorem~2.1 in
  \cite{HSS} shows that $L$, $E$, $E^2$, and $E'$ are all regular with respect
  to $B_1,\dotsc,B_\ell$.  Hence, Proposition~\ref{p:regCrit} together with
  Theorem~2.1 in \cite{HSS} imply that $L$, $E$, and $E'$ satisfy $N_1$, that
  $E'$ has a linear free presentation with respect to $E$, and that both $E$ and
  $E^2$ have a linear free presentation with respect to $E'$.  Therefore,
  Theorem~\ref{t:detPres} proves that $L$ is determinantally presented.
\end{proof}

Theorem~\ref{t:detPres}, combined with results from \S\ref{s:linsyz}, also
yields a proof for our main theorem.

\begin{proof}[Proof of Theorem~\ref{t:main}]
  Let $X$ be a connected scheme of dimension $n$ and let $B$ be a globally
  generated ample line bundle on $X$.  Choose a positive integer $m \in \NN$
  such that $B^m$ is regular with respect to $B$.  Lemma~\ref{l:lower} implies
  that there exists a line bundle $E$, which we may assume is very ample, such
  that, for any nef line bundle $C$, $E \otimes C$ is regular with respect to
  $B$.  By replacing $E$ with $E \otimes B$ if necessary, we may assume that the
  map $\Gamma(B) \otimes H^0(X,E \otimes B^{-1}) \rTo H^0(X,E)$ is surjective.
  Since a sufficiently ample line bundle on $X$ satisfies $N_1$ (combine
  Lemmata~1.1--1.3 in \cite{inamdar} with Fujita's Vanishing Theorem), we may
  also assume that, for any nef line bundle $C$, $E \otimes C$ satisfies $N_1$.

  Consider the line bundle $A := E \otimes B^{2(m-1)n+1}$.  If $L$ is a line
  bundle on $X$ such that $L \otimes A^{-1}$ is nef, then $L = A \otimes C = (E
  \otimes C) \otimes B^{2(m-1)n+1}$ for some nef line bundle $C$.  Our choice of
  $E$ guarantees that, for all $j \geq 1$, $(E \otimes C)^j$ is regular with
  respect to $B$, and that $L$ satisfies $N_1$.  Hence, Proposition~\ref{p:lfp}
  implies that $B^{2(m-1)n+1}$ has a linear free presentation with respect to $E
  \otimes C$.  Proposition~\ref{p:regCrit} together with Mumford's Lemma (e.g.{}
  Theorem~1.8.5 in \cite{PAG}) imply that both $E \otimes C$ and $(E \otimes
  C)^2$ have a linear free presentation with respect to $B^{2(m-1)n+1}$.  Via
  Lemma~\ref{l:regPres} and Proposition~\ref{p:regCrit}, $B^{2(m-1)n+1}$ 
  satisfies $N_1.$ 
  Therefore, Theorem~\ref{t:detPres} proves that $L$ is determinantally
  presented.
\end{proof}

\section{Effective Bounds}
\label{s:effective}

\noindent
In this section, we give effective bounds for determinantally presented line
bundles.  As a basic philosophy, one can convert explicit conditions for line
bundles to satisfy $N_2$ into effective descriptions for determinantally
presented line bundles.  The three subsections demonstrate this philosophy for
products of projective spaces, projective Gorenstein toric varieties, and smooth
varieties.  Despite not being developed here, we expect similar results for
general surfaces and abelian varieties following \cite{GP} and \cites{Rubei,PP}
respectively.

\subsection{Products of Projective Space}
\label{s:proj}

The tools from \S\ref{s:minors} lead to a description of the determinantally
presented ample line bundles on a product of projective spaces.  In contrast
with Theorem~3.11 in \cite{Bernardi} which proves that Segre-Veronese varieties
are defined by $2$\nobreakdash-minors of an appropriate hypermatrix, our
classification shows that a Segre-Veronese variety is typically generated by the
$2$-minors of a single matrix.  In particular, we recover the Segre-Veronese
ideals considered in \cite{Sullivant}*{\S6.2}.

To study the product of projective spaces $X = \PP^{n_1} \times \dotsb \times
\PP^{n_\ell}$, we first introduce some notation.  Let $R := \kk[x_{i,j} : 1 \leq
i \leq \ell, 0 \leq j \leq n_i]$ be the total coordinate ring (a.k.a.{} Cox
ring) of $X$; this polynomial ring has the $\ZZ^\ell$-grading induced by
$\deg(x_{i,j}) := \bm{e}_i \in \ZZ^\ell$.  Hence, we have $R_{\bm{d}} =
\Gamma\bigl( \sO_X(\bm{d}) \bigr)$ for all $\bm{d} \in \ZZ^\ell$ and a
torus-invariant global section of $\sO_X(\bm{d})$ is identified with a monomial
$\bm{x}^{\bm{w}} \in R_{\bm{d}}$ where $\bm{w} \in \NN^r$ and $r :=
\sum_{i=1}^{\ell} (n_i +1)$.  We write $\bm{e}_{i,j}$ for the standard basis of
$\ZZ^r$; in particular $\bm{x}^{\bm{e}_{i,j}} = x_{i,j}$.

\begin{theorem}
  \label{t:proj}
  Let $X = \PP^{n_1} \times \dotsb \times \PP^{n_\ell}$.  An ample line bundle
  $\sO_X(\bm{m})$ is determinantally presented if at least $\ell-2$ of the
  entries in the vector $\bm{m}$ are at least $2$.
\end{theorem}

When $\ell = 2$, this theorem shows that all of the Segre-Veronese embeddings
are determinantally presented.  We note that Corollary~\ref{c:eff} establishes
that $\sO_X(\bm{m})$ is determinantally presented when $m_j \geq 2$ for all $1
\leq j \leq \ell$.

\begin{proof}
  Since a line bundle $\sO_X(\bm{v})$ is ample (and very ample) if and only if
  $v_j \geq 1$ for all $1 \leq j \leq \ell$, Corollary~1.5 in \cite{HSS} shows
  that $\sO_X(\bm{m})$ satisfies $N_1$.  Without loss of generality, we may
  assume that $m_j \geq 2$ for $1 \leq j \leq \ell -2$.  Factor $\sO_X(\bm{m})$
  as $\sO_X(\bm{m}) = E \otimes E'$ where $\bm{u} := \bm{e}_1 + \bm{e}_2 +
  \dotsb + \bm{e}_{\ell-1} = (1,1,\dotsc,1,0)$, $E := \sO_X(\bm{u})$, and $E' :=
  \sO_X(\bm{m}-\bm{u})$.  The canonical surjection $\Gamma(E) \otimes \Gamma(E')
  \rTo \Gamma \bigl( \sO_X(\bm{m}) \bigr)$ implies that the map $\varphi_{\, 2}$
  in \eqref{e:diagram} is surjective.  By Lemma~\ref{l:diagram}, it suffices
  prove that the map $\psi \colon Q_2 \rTo Q_1$ is surjective.  A slight
  modification to the proof of Lemma~4.1 in \cite{Sturmfels} shows that $Q_1 =
  \Ker(\varphi_1)$ is generated by `binomial' elements in $\Sym_2 \bigl(
  \Gamma(E) \bigr) \otimes \Sym_2 \bigl( \Gamma(E') \bigr)$ of the form
  $\bm{x}^{\bm{a}} \bm{x}^{\bm{b}} \otimes \bm{x}^{\bm{c}} \bm{x}^{\bm{d}} -
  \bm{x}^{\bm{a}'} \bm{x}^{\bm{b}'} \otimes \bm{x}^{\bm{c}'} \bm{x}^{\bm{d}'}$
  where $\bm{x}^{\bm{a}}, \bm{x}^{\bm{b}}, \bm{x}^{\bm{a}'}, \bm{x}^{\bm{b}'}
  \in \Gamma(E)$, $\bm{x}^{\bm{c}}, \bm{x}^{\bm{d}}, \bm{x}^{\bm{c}'},
  \bm{x}^{\bm{d}'} \in \Gamma(E')$, and $\bm{a} + \bm{b} + \bm{c} + \bm{d} =
  \bm{a}' + \bm{b}' + \bm{c}' + \bm{d}'$.  Thus, the two terms in each such
  binomial differ by exchanging variables among the various factors.  Since
  every such binomial element is the sum of binomials that each exchange a
  single pair of variables, it suffices to consider the following two cases.

  In the first case, the pair of variables are exchanged between a section of
  $E$ and a section $E'$.  In particular, there exists some $1 \leq k \leq \ell
  -1$ such that the binomial element has the form $\bm{x}^{\bm{a}}
  \bm{x}^{\bm{b}} \otimes \bm{x}^{\bm{c}} \bm{x}^{\bm{d}} - \bm{x}^{\bm{a} -
    \bm{e}_{k, \alpha} + \bm{e}_{k,\gamma}} \bm{x}^{\bm{b}} \otimes
  \bm{x}^{\bm{c} + \bm{e}_{k,\alpha} - \bm{e}_{k, \gamma}} \bm{x}^{\bm{d}}$
  where $\bm{a} - \bm{e}_{k, \alpha}$ and $\bm{c} - \bm{e}_{k, \gamma}$ are
  nonnegative.  This element is the image of $(\bm{x}^{\bm{a}} \otimes
  \bm{x}^{\bm{c}})( \bm{x}^{\bm{b}} \otimes \bm{x}^{\bm{d}}) - (\bm{x}^{\bm{a} -
    \bm{e}_{k, \alpha} + \bm{e}_{k,\gamma}} \otimes \bm{x}^{\bm{c} +
    \bm{e}_{k,\alpha} - \bm{e}_{k, \gamma}} )( \bm{x}^{\bm{b}} \otimes
  \bm{x}^{\bm{d}})$ which lies in $Q_2 = \Ker(\varphi_{\, 2}) \subset \Sym_2
  \bigl( \Gamma(E) \otimes \Gamma(E') \bigr)$.

  In the second case, we may assume that the pair of variables are exchanged
  between two sections of $E'$, as exchanging variables between two sections of
  $E$ is analogous.  More precisely, let $x_{k,\gamma}$ and $x_{k,\delta}$ for
  some $1 \leq k \leq \ell$ denote the exchanged variables and consider the
  binomial element $\bm{x}^{\bm{a}} \bm{x}^{\bm{b}} \otimes \bm{x}^{\bm{c}}
  \bm{x}^{\bm{d}} - \bm{x}^{\bm{a}} \bm{x}^{\bm{b}} \otimes \bm{x}^{\bm{c} -
    \bm{e}_{k, \gamma} + \bm{e}_{k,\delta}} \bm{x}^{\bm{d} + \bm{e}_{k,\gamma} -
    \bm{e}_{k, \delta}}$ where $\bm{c} - \bm{e}_{k, \gamma}$ and $\bm{d} -
  \bm{e}_{k, \delta}$ are nonnegative.  Since $\bm{x}^{\bm{a}} \bm{x}^{\bm{b}}
  \otimes \bm{x}^{\bm{c}} \bm{x}^{\bm{d}} = \bm{x}^{\bm{a}} \bm{x}^{\bm{b}}
  \otimes \bm{x}^{\bm{d}} \bm{x}^{\bm{c}}$ in $\Sym_2\bigl( \Gamma(E) \bigr)
  \otimes \Sym_2 \bigl( \Gamma(E') \bigr)$, we may also assume that $k < \ell$.
  Hence, there is a variable $x_{k,\alpha}$ such that $\bm{a} -
  \bm{e}_{k,\alpha}$ is nonnegative and
  \begin{multline*}
    \bm{x}^{\bm{a}} \bm{x}^{\bm{b}} \otimes \bm{x}^{\bm{c}} \bm{x}^{\bm{d}} -
    \bm{x}^{\bm{a}} \bm{x}^{\bm{b}} \otimes \bm{x}^{\bm{c} - \bm{e}_{k, \gamma}
      + \bm{e}_{k,\delta}} \bm{x}^{\bm{d} + \bm{e}_{k,\gamma} - \bm{e}_{k, \delta}} \\
    \begin{aligned}
      & = \bm{x}^{\bm{a}} \bm{x}^{\bm{b}} \otimes \bm{x}^{\bm{c}}
      \bm{x}^{\bm{d}} - \bm{x}^{\bm{a} - \bm{e}_{k,\alpha} + \bm{e}_{k,\delta}}
      \bm{x}^{\bm{b}} \otimes \bm{x}^{\bm{c}} \bm{x}^{\bm{d} + \bm{e}_{k,\alpha}
        - \bm{e}_{k, \delta}} \\
      & \relphantom{=} + \bm{x}^{\bm{a} - \bm{e}_{k,\alpha} + \bm{e}_{k,\delta}}
      \bm{x}^{\bm{b}} \otimes \bm{x}^{\bm{c}} \bm{x}^{\bm{d} + \bm{e}_{k,\alpha}
        - \bm{e}_{k, \delta}} - \bm{x}^{\bm{a} - \bm{e}_{k,\alpha} +
        \bm{e}_{k,\gamma}} \bm{x}^{\bm{b}} \otimes \bm{x}^{\bm{c} -
        \bm{e}_{k,\gamma} + \bm{e}_{k, \delta}} \bm{x}^{\bm{d} +
        \bm{e}_{k,\alpha} - \bm{e}_{k, \delta}} \\
      & \relphantom{=} + \bm{x}^{\bm{a} - \bm{e}_{k,\alpha} + \bm{e}_{k,\gamma}}
      \bm{x}^{\bm{b}} \otimes \bm{x}^{\bm{c} - \bm{e}_{k,\gamma} + \bm{e}_{k,
          \delta}} \bm{x}^{\bm{d} + \bm{e}_{k,\alpha} - \bm{e}_{k, \delta}} -
      \bm{x}^{\bm{a}} \bm{x}^{\bm{b}} \otimes \bm{x}^{\bm{c} - \bm{e}_{k,\gamma}
        + \bm{e}_{k, \delta}} \bm{x}^{\bm{d} + \bm{e}_{k,\gamma} - \bm{e}_{k,
          \delta}} \, .
    \end{aligned}
  \end{multline*}
  In other words, the binomial element under consideration is a sum of binomials
  in which variables are exchanged between sections of $E$ and $E'$.  Hence, the
  first case shows that this binomial element lies in the image of $Q_2$.

  We conclude that $\psi$ is surjective and $\sO_X(\bm{m})$ is determinantally
  presented.
\end{proof}

The next proposition shows that Theorem~\ref{t:proj} is optimal when $\ell = 3$.
In fact, our experiments in \emph{Macaulay2}~\cite{M2} suggest that
Theorem~\ref{t:proj} is always sharp.

\begin{proposition}
  If $X = \PP^{n_1} \times \dotsb \times \PP^{n_\ell}$ with $\ell \geq 3$, then
  the ample line bundle $\sO_X(\bm{1})$ is not determinantally presented.
\end{proposition}

\begin{proof}
  Any nontrivial factorization of $\sO_X(\bm{1})$ has the form $E \otimes E'$
  where $E := \sO_X(\bm{u})$ for some $\bm{u} \in \{0,1\}^\ell$ and $E' :=
  \sO_X(\bm{1} - \bm{u})$.  For a suitable choice of bases for $\Gamma\bigl(
  \sO_X(\bm{u})\bigr)$ and $\Gamma\bigl( \sO_X(\bm{1} - \bm{u})\bigr)$, the
  associated matrix $\Omega\bigl( \sO_X(\bm{u}), \sO_X(\bm{1}-\bm{u}) \bigr)$ is
  the generic $(s \times t)$-matrix with $s := \sum_{u_i \neq 0}(n_i+1)$ and $t
  := \sum_{i=0}^{\ell} (n_i+1) - s$.  Since the $2$-minors of a generic $(s
  \times t)$-matrix define $\PP^{s-1} \times \PP^{t-1}$ in its Segre embedding,
  we see that $\sO_X(\bm{1})$ is not determinantally presented when $\ell \geq
  3$.
\end{proof}

\begin{example}
  \label{e:PP1^3}
  Consider the variety $X = \PP^1 \times \PP^1 \times \PP^1$ embedded in
  $\PP^{11} = \Proj(\kk[y_0,\dotsc,y_{11}])$ by the complete linear series of
  $\sO_X(2,1,1)$. If $R = \kk[x_{1,0}, x_{1,1}, x_{2,0}, x_{2,1}, x_{3,0},
  x_{3,1}]$ is the total coordinate ring of $X$, then the twelve monomials
  \[
  \left\{ \begin{matrix} x_{1,0}^2x_{2,0}x_{3,0}, & x_{1,0}^2x_{2,0}x_{3,1},
      & x_{1,0}^2x_{2,1}x_{3,0}, & x_{1,0}^2x_{2,1}x_{3,1}, \\
      x_{1,0}x_{1,1}x_{2,0}x_{3,0}, & x_{1,0}x_{1,1}x_{2,0}x_{3,1}, &
      x_{1,0}x_{1,1}x_{2,1}x_{3,0}, & x_{1,0}x_{1,1}x_{2,1}x_{3,1}, \\
      x_{1,1}^2x_{2,0}x_{3,0}, & x_{1,1}^2x_{2,0}x_{3,1}, &
      x_{1,1}^2x_{2,1}x_{3,0}, & x_{1,1}^2x_{2,1}x_{3,1}
    \end{matrix} \right\}
  \]
  give an ordered basis for $\Gamma\bigl( \sO_X(2,1,1) \bigr)$. The homogeneous
  ideal $I_{X | \PP^{11}}$ is the toric ideal associated to these monomials and
  is minimally generated by thirty three quadrics.  Choosing $\{ x_{1, 0}x_{2,
    0}, x_{1, 0}x_{2, 1}, x_{1, 1}x_{2, 0}, x_{1, 1}x_{2, 1} \}$ and $\{ x_{1,
    0}x_{3, 0}, x_{1, 0}x_{3, 1}, x_{1, 1}x_{3, 0}, x_{1, 1}x_{3, 1} \}$ as
  ordered bases for $\Gamma \bigl( \sO_X(1,1,0) \bigr)$ and $\Gamma \bigl(
  \sO_X(1,0,1) \bigr)$, $\Omega\bigl( \sO_X(1,1,0), \sO_X(1,0,1) \bigr)$ is
  \[
  \left[ \begin{matrix}
      y_0 & y_1 & y_4 & y_5 \\
      y_2 & y_3 & y_6 & y_7 \\
      y_4 & y_5 & y_8 & y_9 \\
      y_6 & y_7 & y_{10} & y_{11}
    \end{matrix} \right]
  \]
  and one may verify that the 2-minors of this matrix generates the ideal of
  $X$, so $\sO_X(2,1,1)$ is determinantally presented. \hfill $\diamond$
\end{example}

However, if we consider multiple factorizations of a very ample line bundle on a
product of projective spaces, then we do obtain a convenient expression of the
homogeneous ideal as the $2$-minors of matrices.  This perspective give a
conceptual explanation for both Theorem~2.6 in \cite{Ha} and Theorem~3.11 in
\cite{Bernardi}.

\begin{proposition}
  \label{p:mult}
  If $X = \PP^{n_1} \times \dotsb \times \PP^{n_\ell}$, then the homogeneous
  ideal of $X$ in $\PP\bigl( \sO_X(\bm{d}) \bigr)$ is generated by the
  $2$-minors of the matrices $\Omega\bigl( \sO_X(\bm{e}_i),
  \sO_X(\bm{d}-\bm{e}_i) \bigr)$ where $1 \leq i \leq \ell$.
\end{proposition}

\begin{proof}
  Given Theorem~\ref{t:proj}, we may assume that $\ell \geq 3$.  For brevity,
  set $E_i := \sO_X(\bm{e}_i)$ and $E'_i := \sO_X(\bm{d} - \bm{e}_i)$ where $1
  \leq i \leq \ell$.  Since $\Gamma(E_i) \otimes \Gamma(E'_i)$ surjects onto
  $\Gamma \bigl( \sO_X(\bm{d}) \bigr)$, the map $\varphi_{\, 2}$ in
  \eqref{e:diagram} is surjective, and it suffices prove that the map $\psi
  \colon Q_2 \rTo Q_1$ is surjective. By an abuse of notation, we use
  $\epsilon_i$ to denote the canonical inclusion map onto the $i$-th summand for
  all three of the direct sums appearing in the middle column of
  \eqref{e:diagram}.  As in the proof of Theorem~\ref{t:proj}, $Q_1$ is
  generated by binomial elements in $\bigoplus_{k=1}^{\ell} \Sym_2 \bigl(
  \Gamma(E_k) \bigr) \otimes \Sym_2 \bigl( \Gamma(E_k') \bigr)$.  Generators
  have the form $\epsilon_i(x_{i,\alpha} x_{i,\beta} \otimes \bm{x}^{\bm{c}}
  \bm{x}^{\bm{d}}) - \epsilon_j(x_{j,\gamma} x_{j,\delta} \otimes
  \bm{x}^{\bm{a}} \bm{x}^{\bm{b}})$ where $x_{i,\alpha}, x_{i,\beta} \in
  \Gamma(E_i)$, $\bm{x}^{\bm{c}}, \bm{x}^{\bm{d}} \in \Gamma(E'_i)$,
  $x_{j,\gamma}, x_{j,\delta} \in \Gamma(E_j)$, $\bm{x}^{\bm{a}},
  \bm{x}^{\bm{b}} \in \Gamma(E'_j)$ and $\bm{e}_{i,\alpha} + \bm{e}_{i,\beta} +
  \bm{c} + \bm{d} = \bm{a} + \bm{b} + \bm{e}_{j,\gamma} + \bm{e}_{j,\delta}$.
  We consider the following two cases.

  In the first case, we have $i = j$.  Since every binomial element is the sum
  of binomials that each exchange a single pair of variables, it suffices to
  consider an element of the form $\epsilon_i \bigl( x_{i,\alpha} x_{i,\beta}
  \otimes \bm{x}^{\bm{c}} \bm{x}^{\bm{d}} - x_{i,\alpha} x_{i,\beta} \otimes
  \bm{x}^{\bm{c} - \bm{e}_{k,\gamma} + \bm{e}_{k,\delta}} \bm{x}^{\bm{d} +
    \bm{e}_{k,\gamma} - \bm{e}_{k,\delta}} \bigr)$ where $1 \leq k \leq \ell$
  and both $\bm{c} - \bm{e}_{k,\gamma}$ and $\bm{d} - \bm{e}_{k,\delta}$ are
  nonnegative.  This element is the image of
  \begin{multline*}
    \epsilon_i \bigl( (x_{i,\alpha} \otimes \bm{x}^{\bm{c}}) (x_{i,\beta}
    \otimes \bm{x}^{\bm{d}}) - (x_{i,\alpha} \otimes \bm{x}^{\bm{c} -
      \bm{e}_{k,\gamma} + \bm{e}_{k,\delta}}) (x_{i,\beta} \otimes
    \bm{x}^{\bm{d} + \bm{e}_{k,\gamma} -
      \bm{e}_{k,\delta}}) \bigr) \\
    - \epsilon_k \bigl( (x_{k,\gamma} \otimes \bm{x}^{\bm{c} + \bm{e}_{i,\alpha}
      - \bm{e}_{k,\gamma}}) (x_{k,\delta} \otimes \bm{x}^{\bm{d} +
      \bm{e}_{i,\beta} - \bm{e}_{k,\delta}}) - (x_{k,\delta} \otimes
    \bm{x}^{\bm{c} + \bm{e}_{i,\alpha} - \bm{e}_{k,\gamma}}) (x_{k,\gamma}
    \otimes \bm{x}^{\bm{d} + \bm{e}_{i,\beta} - \bm{e}_{k,\delta}}) \bigr) \, ,
  \end{multline*}
  which lies in $Q_2 = \Ker(\varphi_{\, 2})$.
 
  For the second case, we have $i \neq j$.  We may assume that the binomial
  element has the form $\epsilon_i ( x_{i,\alpha} x_{i, \beta} \otimes
  \bm{x}^{\bm{c}} \bm{x}^{\bm{d}} ) - \epsilon_j ( x_{j,\gamma} x_{j, \delta}
  \otimes \bm{x}^{\bm{c} + \bm{e}_{i, \alpha} - \bm{e}_{j,\gamma}}
  \bm{x}^{\bm{d} + \bm{e}_{i,\beta} - \bm{e}_{j,\delta}})$, where $\bm{c}-
  \bm{e}_{j,\gamma}$ and $\bm{d} - \bm{e}_{j,\delta}$ are nonnegative, because
  any additional exchanges of variables can be obtained by adding elements from
  the first case.  This element is the image of
  \[
  \epsilon_i \bigl( (x_{i,\alpha} \otimes \bm{x}^{\bm{c}}) (x_{i,\beta} \otimes
  \bm{x}^{\bm{d}}) \bigr) - \epsilon_j \bigl( (x_{j,\gamma} \otimes
  \bm{x}^{\bm{c} + \bm{e}_{i,\alpha} - \bm{e}_{j,\gamma}}) (x_{j,\delta} \otimes
  \bm{x}^{\bm{d} + \bm{e}_{i,\beta} - \bm{e}_{j,\delta}}) \bigr) \, ,
  \]
  which lies in $Q_2 = \Ker(\varphi_{\, 2})$.
\end{proof}

\begin{example}
  Consider the variety $X = \PP^1 \times \PP^1 \times \PP^1$ embedded in $\PP^7
  = \Proj(\kk[y_0,\dotsc,y_7])$ by the complete linear series of the line bundle
  $\sO_X(1,1,1)$. The homogeneous ideal $I_{X | \PP^7}$ is the toric ideal
  associated to the monomial list
  \[
  \left\{ 
    \begin{matrix}
      x_{1,0}x_{2,0}x_{3,0}, & x_{1,0}x_{2,0}x_{3,1}, &
      x_{1,0}x_{2,1}x_{3,0}, & x_{1,0}x_{2,1}x_{3,1}, \\
      x_{1,1}x_{2,0}x_{3,0}, & x_{1,1}x_{2,0}x_{3,1}, &
      x_{1,1}x_{2,1}x_{3,0}, & x_{1,1}x_{2,1}x_{3,1}
    \end{matrix}
  \right\}
  \]
  and is minimally generated by nine quadrics.  Choosing appropriate monomials
  for the ordered bases of the global sections, we obtain
  \begin{align*}
    \Omega\bigl( \sO_X(1,0,0), \sO_X(0,1,1) \bigr) &= \left[ 
      \begin{matrix}
        y_0 & y_1 & y_2 & y_3 \\
        y_4 & y_5 & y_6 & y_7 
      \end{matrix} \right], \\
    \Omega\bigl( \sO_X(0,1,0), \sO_X(1,0,1) \bigr) &= \left[ 
      \begin{matrix}
        y_0 & y_1 & y_4 & y_5 \\
        y_2 & y_3 & y_6 & y_7
      \end{matrix} \right], \\
    \Omega\bigl( \sO_X(0,0,1), \sO_X(1,1,0) \bigr) &= \left[ 
      \begin{matrix}
        y_0 & y_2 & y_4 & y_6 \\
        y_1 & y_3 & y_5 & y_7 
      \end{matrix} \right] \, .
  \end{align*}
  It follows that $\sO_X(1,1,1)$ is not determinantally presented, but one
  easily verifies that the ideal $I_{X | \PP^7}$ is generated by the $2$-minors
  of all three matrices. \hfill $\diamond$
\end{example}

Multiple factorizations of a very ample line bundle allow one to describe a
larger number of homogeneous ideals via $2$\nobreakdash-minors.  With this in
mind, it would be interesting to write down the analogue of
Theorem~\ref{t:detPres} for multiple factorizations of the line bundle.

\subsection{Toric Varieties}

In addition to the bound given in Corollary~\ref{c:eff}, there is an effective
bound for toric varieties involving adjoint bundles for toric varieties; cf.{}
Corollary~1.6 in \cite{HSS}.  Recall that a line bundle on a toric variety $X$
is nef if and only if it is globally generated, and the dualizing sheaf $K_X$ is
a line bundle if and only if $X$ is Gorenstein.

\begin{proposition}
  \label{p:toric}
  Let $X$ be a projective $n$-dimensional Gorenstein toric variety with
  dualizing sheaf $K_X$, and let $B_1, \dotsc, B_\ell$ be the minimal generators
  of its nef cone $\Nef(X)$.  Suppose that $\bm{m}, \bm{m}' \in \NN^\ell$
  satisfy $\bm{B}^{\bm{m} - \bm{u}}, \bm{B}^{\bm{m}' - \bm{u}} \in \fB$ for all
  $\bm{u} \in \NN^\ell$ with $|\bm{u}| \leq n+1$.  If $X \neq \PP^n$ and $\bm{w}
  \in \NN^\ell$, then $L = K_X^2 \otimes \bm{B}^{\bm{m} + \bm{m}' +\bm{w}}$ is
  determinantally presented.
\end{proposition}

\begin{proof}
  Factor $L$ as $L = E \otimes E'$ where $E := K_X \otimes
  \bm{B}^{\bm{m}+\bm{w}}$ and $E' := K_X \otimes \bm{B}^{\bm{m}'}$.  Since
  $\bm{B}^{\bm{m} - (n+1)\bm{e}_j}, \bm{B}^{\bm{m}' - (n+1)\bm{e}_j} \in \fB$,
  Corollary~0.2 in \cite{Fujino} implies that $E \otimes \bm{B}^{-\bm{e}_j}$ and
  $E' \otimes \bm{B}^{-\bm{e}_j}$ belong to $\fB$ for all $1 \leq j \leq \ell$.
  For any torus invariant curve $Y$, there is a $\bm{B}^{\bm{e}_j}$ such that
  $\bm{B}^{\bm{e}_{j}} \cdot Y > 0$.  Theorem~3.4 in \cite{Mustata} implies that
  $E$, $E^2$ and $E'$ are regular with respect to $B_1, \dotsc, B_\ell$.  Hence,
  Proposition~\ref{p:regCrit} shows that $L$, $E$, and $E'$ satisfy $N_1$, and
  that $E$ has a linear free presentation with respect to $E'$, $E'$ has a
  linear free presentation with respect to $E$, and $E^2$ has a linear free
  presentation with respect to $E'$.  Therefore, Theorem~\ref{t:detPres} shows
  that $L$ is determinantally presented.
\end{proof}

\begin{proof}[Proof of Theorem~\ref{t:effective} for toric varieties]
  This is a special case of Proposition~\ref{p:toric}.
\end{proof}

We give an example showing that Theorem~\ref{t:effective} is not sharp for all
toric varieties.

\begin{example}
  Consider the toric del Pezzo surface $X$ obtained by blowing up $\PP^2$ at the
  three torus-fixed points.  Let $R := \kk[x_0, \dotsc, x_5]$ be the total
  coordinate ring of $X$.  The anticanonical bundle $K_X^{-1}$ is very ample
  and corresponds to polygon 
  \[
  P := \conv \bigl\{ (1,0), (1,1), (0,1), (-1,0), (-1,-1), (0,-1) \bigr\} \, .
  \] 
  It is easy to see that the polygon $P$ is the smallest lattice polygon with
  its inner normal fan.  The polygon $2P$ contains 19 lattice points.  The
  corresponding monomials
  \[
  \left\{
    \begin{matrix}
      x_0^4x_1^4x_2^2x_5^2, & x_0^4x_1^3x_2x_4x_5^3, & x_0^4x_1^2x_4^2x_5^4, &
      x_0^3x_1^4x_2^3x_3x_5, & x_0^3x_1^3x_2^2x_3x_4x_5^2, \\
      x_0^3x_1^2x_2x_3x_4^2x_5^3, & x_0^3x_1x_3x_4^3x_5^4, &
      x_0^2x_1^4x_2^4x_3^2, & x_0^2x_1^3x_2^3x_3^2x_4x_5, &
      x_0^2x_1^2x_2^2x_3^2x_4^2x_5^2, \\
      x_0^2x_1x_2x_3^2x_4^3x_5^3, & x_0^2x_3^2x_4^4x_5^4, &
      x_0x_1^3x_2^4x_3^3x_4, & x_0x_1^2x_2^3x_3^3x_4^2x_5, &
      x_0x_1x_2^2x_3^3x_4^3x_5^2, \\
      x_0x_2x_3^3x_4^4x_5^3, & x_1^2x_2^4x_3^4x_4^2, & x_1x_2^3x_3^4x_4^3x_5, &
      x_2^2x_3^4x_4^4x_5^2
    \end{matrix}
  \right\}
  \]
  embed $X$ into $\PP^{18} = \Proj\bigl( \kk[y_0, \dotsc, y_{18}] \bigr)$.  The
  homogeneous ideal $I_{X | \PP^{18}}$ is the toric ideal associated to these
  monomials and is minimally generated by $129$ quadrics.  Choosing $\{
  x_0^2x_1^2x_2x_5, x_0^2x_1x_4x_5^2, x_0x_1^2x_2^2x_3, x_0x_1x_2x_3x_4x_5,
  x_0x_3x_4^2x_5^2, x_1x_2^2x_3^2x_4, x_2x_3^2x_4^2x_5 \}$ as an ordered basis
  for $\Gamma(K_X^{-1}),$ the matrix $\Omega\bigl( K_X^{-1}, K_X^{-1} \bigr)$ is
  \[
  \left[ 
    \begin{matrix}
      y_0 & y_1 & y_3 & y_4 & y_5 & y_8 & y_9 \\
      y_1 & y_2 & y_4 & y_5 & y_6 & y_9 & y_{10} \\
      y_3 & y_4 & y_7 & y_8 & y_9 & y_{12} & y_{13} \\
      y_4 & y_5 & y_8 & y_9 & y_{10} & y_{13} & y_{14} \\
      y_5 & y_6 & y_9 & y_{10} & y_{11} & y_{14} & y_{15} \\
      y_8 & y_9 & y_{12} & y_{13} & y_{14} & y_{16} & y_{17} \\ 
      y_9 & y_{10} & y_{13} & y_{14} & y_{15} & y_{17} & y_{18}
    \end{matrix}
  \right] \, , 
  \]
  and its $2$-minors generate $I_{X |\PP^{18}}.$ However,
  Theorem~\ref{t:effective} only establishes that the line bundle $K_X^{-4} =
  K_X^2 \otimes (K_X^{-1})^{2 \cdot 2 +2}$ is determinantally presented.  \hfill
  $\diamond$
\end{example}

\subsection{Smooth Varieties}

For smooth varieties, we also have an effective bound for adjoint bundles; see
Theorem~\ref{t:effective}.

\begin{proof}[Proof of Theorem~\ref{t:effective} for smooth varieties]
  Factor the line bundle $L$ as $L = E \otimes E'$ where $E := K_X \otimes
  A^{n+1}$ and $E' := K_X \otimes A^{j-n-1} \otimes B$. Since $j \geq 2n+2$ and
  $E$ is nef (see Example~1.5.35 in \cite{PAG}), Proposition~\ref{p:smoothAdj}
  implies that $L$, $E$, and $E'$ satisfy $N_1$, $E$ has a linear free
  presentation with respect to $E'$, $E'$ has a linear free presentation with
  respect to $E$, and $E^2$ has a linear free presentation with respect to
  $E'$. Therefore, Theorem~\ref{t:detPres} shows that $L$ is determinantally
  presented.
\end{proof}

We end with an example showing that the hypotheses in Theorem~\ref{t:effective}
are optimal without further restrictions on the varieties under consideration.

\begin{example}
  Consider the Grassmannian $X = \operatorname{Gr}(2,4)$ parametrizing all two
  dimensional subspaces of the vector space $\kk^4$.  Let $\sO_X(1)$ denote the
  determinant of the universal rank $2$ sub-bundle on $X$.  The associated
  complete linear series determines the Pl\"ucker embedding of $X$ into $\PP^5 =
  \Proj(\kk[x_{1,2},x_{1,3},x_{1,4},x_{2,3},x_{2,4},x_{3,4}])$.  Since
  $I_{X|\PP^5} = \langle x_{1,2}x_{3,4} - x_{1,3}x_{2,4} + x_{2,3}x_{1,4}
  \rangle$, it follows that $\sO_X(1)$ is not determinantally presented.  On the
  other hand, the monomials
  \[
  \left\{ 
    \begin{matrix} 
      x_{1,2}^2, & x_{1,2}x_{1,3}, & x_{1,2}x_{1,4}, &
      x_{1,2}x_{2,3}, & x_{1,2}x_{2,4}, & x_{1,2}x_{3,4}, & x_{1,3}^2, \\
      x_{1,3}x_{1,4}, & x_{1,3}x_{2,3}, & & x_{1,3}x_{3,4}, & x_{1,4}^2, &
      x_{1,4}x_{2,3}, & x_{1,4}x_{2,4}, \\ x_{1,4}x_{3,4}, & x_{2,3}^2, &
      x_{2,3}x_{2,4}, & x_{2,3}x_{3,4}, & x_{2,4}^2, & x_{2,4}x_{3,4}, &
      x_{3,4}^2
    \end{matrix}
  \right\}
  \]
  form an ordered basis for $\Gamma\bigl( \sO_X(2) \bigr)$, so the complete
  linear series of $\sO_X(2)$ embeds $X$ into $\PP^{19} = \Proj(\kk[y_0,\dotsc,
  y_{19}])$.  The matrix $\Omega \bigl( \sO_X(1), \sO_X(1) \bigr)$ is
  \[
  \begin{bmatrix}
    y_0 & y_1 & y_2 & y_3 & y_4 & y_5 \\
    y_1 & y_6 & y_7 & y_8 & y_5+y_{11} & y_9 \\
    y_2 & y_7 & y_{10} & y_{11} & y_{12} & y_{13} \\
    y_3 & y_8 & y_{11} & y_{14} & y_{15} & y_{16} \\
    y_4 & y_5+y_{11} & y_{12} & y_{15} & y_{17} & y_{18} \\ 
    y_5 & y_9 & y_{13} & y_{16} & y_{18} & y_{19} 
  \end{bmatrix} 
  \]
  and the $2$-minors of this matrix generated $I_{X | \PP^{19}}$ (indeed, this
  is the second Veronese of the Pl\"ucker embedding).  Since $K_X = \sO_X(-4)$
  and $\sO_X(2) = K_X^2 \otimes \sO_X(1)^{2 \cdot 4+2}$, we see that the bound
  in Theorem~\ref{t:effective} is sharp in this case.  \hfill $\diamond$
\end{example}

\raggedright

\end{document}